\documentclass[11pt]{article}

\usepackage[margin=1in]{geometry} 
\usepackage{amsmath,amsthm,amssymb,amsfonts}
\usepackage{graphicx, multicol, array}
\usepackage{cases}
\usepackage{enumerate}
\usepackage{enumitem}
\usepackage{hyperref}

\usepackage{mathrsfs}

\newcommand{\Z}{\mathbb{Z}}
\newcommand{\Q}{\mathbb{Q}}
\newcommand{\R}{\mathbb{R}} 
\newcommand{\C}{\mathbb{C}}
\newcommand{\F}{\mathbb{F}}
\newcommand{\A}{\mathbb{A}}
\newcommand{\T}{\mathbb{T}}
\newcommand{\I}{\mathbb{I}}

\newtheorem{thm}{Theorem}[section]
\newtheorem{lem}{Lemma}[section]

\newtheorem{exa}{Example}[section]
\newtheorem{cor}{Corollary}[section]
\newtheorem{dfn}{Definition}[section]
\newtheorem{rem}{Remark}[section]
\newtheorem{exe}{Exercise}[section]

\title{{\Large Irrationality of the Zeta Constants}}
\date{}
\author{N. A. Carella}

\usepackage{fancyhdr}
\begin{document}
\maketitle

\vskip .25 in 
\textbf{Abstract:} A general technique for proving the irrationality of the zeta constants \(\zeta (2n+1), n \geq  1\), from 	the known irrationality of the beta constants\textit{  }\(L(2n+1,\chi )\) is developed in this note. The irrationality of the zeta constants \(\zeta(2n), n \geq  1\), and \(\zeta(3)\) are well known, but the irrationality results for the zeta constants \(\zeta (2n+1), n \geq 2\), are new, and seem to show that these are irrational numbers. By symmetry, the irrationality of the beta constants \(L(2n,\chi )\) are derived from the known rrationality of the zeta constants 

\(\zeta (2n)\).  
\\


\vskip .25 in 


 \textbf{AMS Mathematical Subjects Classification:} 11J72; 11A55.\\

\vskip .75 in
\section{Introduction}\label{section1}

As the rationality or irrationality nature of a number is an arithmetic property, it is not surprising to encounter important constants, whose rationality or irrationality is linked to the properties of the integers and the distribution of the prime numbers, for example, the number \(6/\pi ^{2}=\prod_{p\geq 2} \left(1-p^{-2}\right)\). An estimate of the partial sum of the Dedekind zeta function of quadratic numbers fields will be utilized to develop a general technique for proving the irrationality of the zeta constants \(\zeta (2n+1)\) from the known irrationality of the beta constants \(L(2n+1,\chi ), 1 \neq  n \in  \mathbb{N}\). This technique provides another proof of the first odd case \(\zeta (3)\), which have well known proofs of irrationalities, see \cite{AR79}, \cite{BF79}, \cite{SJ05}, et al, and an original proof for the other odd cases \(\zeta (2n+1), n \geq  2\), which seems to confirm the irrationality of these number.\\

\begin{thm}\label{thm1.1}
For each fixed odd integer \(s = 2k + 1 \geq  3\), the zeta constant \(\zeta (s)\) is an irrational number.
\end{thm}

The current research literature on the zeta constants \(\zeta (2n+1)\) states the following:
\begin{enumerate}
\item The special zeta value \(\zeta (3)\) is an irrational number, see \cite{AR79}, \cite{BF79}, \cite{HD01}, \cite{SJ05}, et al.

\item At least one of the four numbers \(\zeta (5), \zeta (7), \zeta (9),\) and \(\zeta (11)\) is an irrational number, see \cite[p.\ 7]{WM00} and 

\cite{WZ04}.
\item The sequence \(\zeta (5), \zeta (7), \zeta (9), \zeta (11), \ldots\) contains infinitely many irrational zeta constants, see \cite{BR01}. Various 

advanced techniques for studying the zeta constants are surveyed in \cite{NY09}, and \cite{WM05}. 
\end{enumerate}

By the symmetry of the factorization of the Dedekind zeta function \(\zeta _K(s)=\zeta (s)L(s,\chi )\) with respect to either \(\zeta (s)\) or \(L(s,\chi)\), almost the same analysis leads to a derivation of the irrationality of the beta constants \(L(2n,\chi )\) from the known irrationality of the zeta constants \(\zeta (2n)\) \(n\geq 1\).\\

\begin{thm}\label{thm1.2}
For each fixed even integer \(s=2n\geq 2\), and the nonprincipal quadratic character $\chi \mod 4$, the beta constant \(L(s,\chi)=\sum _{n=1}^{\infty } \chi (n)n^{-s}\) is an irrational number.
\end{thm}

 Section \ref{sec2} contains basic materials on the theory of irrationality, the transcendental properties of real numbers, and an estimate of the summatory function of the Dedekind zeta function. The proof of Theorem \ref{thm1.1} is given in section \ref{sec6}. Last but not least, a sketch of the proof of Theorem \ref{thm1.2} is given in the last Section.\\

\section{Fundamental Concepts and Background}\label{sec2}
The basic notation, concepts and results employed throughout this work are stated in this Section. All the materials covered in this subsection are standard definitions and results in the literature, confer \cite{HW08}, \cite{LS95}, \cite{NZ91}, \cite{SJ05}, \cite{WM00}, et al. 

\subsection{Criteria for Rationality and for Irrationality}
A real number \(\alpha \in \mathbb{R}\) is called \textit{rational} if \(\alpha = a/b\), where \(a, b \in \mathbb{Z}\) are integers. Otherwise, the number
is \textit{ irrational}. The irrational numbers are further classified as \textit{algebraic} if \(\alpha\) is the root of an irreducible polynomial \(f(x) \in \mathbb{Z}[x]\) of degree \(\deg (f)>1\), otherwise it is \textit{transcendental}.\\

\begin{lem} \label{lem2.1} {\normalfont  (Criterion for rationality)} If a real number \(\alpha \in \mathbb{Q}\) is a rational number, then there exists a constant \(c = c(\alpha )\) such that
\begin{equation}
\frac{c}{q}\leq \left|  \alpha -\frac{p}{q} \right|
\end{equation}
holds for any rational fraction \(p/q \neq \alpha\). Specifically, \(c \geq  1/B\text{ if }\alpha = A/B\).
\end{lem}

This is a statement about the lack of effective or good approximations of an arbitrary rational number \(\alpha \in \mathbb{Q}\) by other rational numbers. On the other hand, irrational numbers \(\alpha \in \mathbb{R}-\mathbb{Q}\) have effective approximations by rational numbers.\\

If the complementary inequality \(\left|  \alpha -p/q \right| <c/q\) holds for infinitely many rational approximations \(p/q\), then it already shows that the real number \(\alpha \in \mathbb{R}\) is almost irrational, so it is almost sufficient to prove the irrationality of real numbers.

\begin{lem} \label{lem2.2} {\normalfont  (Criterion for irrationality)} Let $\psi(x)=o(1/x)$ be a monotonic decresing function, and let \(\alpha \in \mathbb{Q}\) be a real number. If
\begin{equation}
0< \left|  \alpha -\frac{p}{q} \right|< \psi(q)
\end{equation}
holds for infinitely many rational fraction \(p/q \in \mathbb{Q}\), then \(\alpha \) is irrational.
\end{lem}

\begin{proof} By Lemma \ref{lem2.1} and the hypothesis, it follows that

\begin{equation}
\frac{c}{q}\leq \left|  \alpha -\frac{p}{q} \right| < \psi(q)=o\left (\frac{1}{q} \right).
\end{equation}
But this is a contradiction since $c/q \ne o(1/q)$.
\end{proof}

More precise results for testing the irrationality of an arbitrary real number are stated below.\\

\begin{thm} \label{thm2.1} {\normalfont (Dirichlet)} Suppose $\alpha \in \mathbb{R}$ is an irrational number. Then there exists an infinite
sequence of rational numbers $p_n/q_n$ satisfying
\begin{equation}
0 < \left|  \alpha -\frac{p_n}{q_n} \right|< \frac{1}{q_n^2}
\end{equation}
for all integers \(n\in \mathbb{N}\).
\end{thm}

For continued fractions $\alpha=[a_0,a_1,a_2, \ldots]$ with sizable entries $a_i \geq a>1$, where $a$ ia a constant, there is a slightly better inequality.\\

\begin{thm} \label{thm2.2}   Let $\alpha=[a_0,a_1,a_2, \ldots]$ be the continued fraction of a real number, and let $\{p_n/q_n: n \geq 1\}$ be the sequence of convergents. Then
\begin{equation}
0 < \left|  \alpha -\frac{p_n}{q_n} \right|< \frac{1}{a_nq_n^2}
\end{equation}
for all integers \(n\in \mathbb{N}\).
\end{thm}
This is standard in the literature, the proof appears in \cite[Theorem 171]{HW08}, \cite[Corollary 3.7]{SJ05}, and similar references. A more general result provides a family of inequalities for almost all real numbers.

\begin{thm} \label{thm2.3} {\normalfont (Khinchin)} Let \(\psi \) be a real and deceasing function, and let $\alpha \in \mathbb{R}$ be a real number. If there exists an infinite sequence of rational
approximations \(p_n/q_n\) such that \(p_n/q_n\neq \alpha\), and
\begin{equation}
\left|  \alpha -\frac{p_n}{q_n} \right|<\frac{\psi(q_n)}{q_n}
\end{equation}
and $\sum_{q} \psi(q) < \infty$, then real number $\alpha$ is $\psi$-approximatable.
\end{thm}

\section{Estimate of An Arithmetic Function }\label{sec3}
Let \(q \geq  2\) be an integer, and let\(\text{  }\chi \neq 1\) be the quadratic character modulo \(q\). The Dedekind zeta function of a quadratic numbers field \(\mathbb{Q}\left(\sqrt{q}\right)\) is defined by \(\zeta _K(s)=\zeta (s)L(s,\chi )\). The factorization consists of the zeta function
{ }\(\zeta (s)=\sum _{n=1}^{\infty } n^{-s}\), and the $L$-function \(L(s,\chi )=\sum _{n=1}^{\infty } \chi (n)n^{-s}\), see \cite[p.\ 219]{CO07}. The product has the Dirichlet series expansion 
\begin{equation}\label{400}
\zeta (s)\cdot L(s,\chi )=\left(\sum _{n=1}^{\infty } \frac{1}{n^s}\right)\left(\sum _{n=1}^{\infty } \chi
(n)\frac{1}{n^s}\right)=\sum _{n=1}^{\infty } \left(\sum _{d | n} \chi (d)\right)\frac{1}{n^s}, 
\end{equation}
and its summatory function is \(\sum _{n\leq x} r(n)=4\sum
_{n\leq x}\sum _{d | n} \chi (d)\), see \cite[p.\ 17]{IK04}.  \\

The finite sum $r_Q(n)=4\sum _{d | n} \chi (d)$ is the number of respresentations of $n$ by primitive forms $Q(u,v)=au^2+buv+cv^2$ of discriminant $q$. For $q=4$, the nonprincipal character $\chi \ne 1$ has order 2, and the counting function \(r(n)=4\sum _{d | n} \chi (d)=\#\left\{ (u,v):n=u^2+v^2 \right\}\geq 0\)
tallies the number of representations of an integer \(n \geq  1\) as sums of two squares.\\

\begin{lem} \label{lem3.1} The average order of the summatory function of the Dedekind function is given by
\begin{equation}
\sum _{n \leq  x} \frac{r_Q(n)}{n^s}=\zeta (s)L(s,\chi )+\frac{c_0}{x^{s-1}}+O\left(\frac{1}{x^{s-1/2}}\right),
\end{equation}
where \(c_0>0\) is a constant, and \(x\geq 1\) is a large number.
\end{lem}

A proof appears in \cite[p.\ 369]{MV07}, slightly different version, based on the hyperbola method, is computed in \cite[p.\ 255]{MR08}.

\section{The Irrationality of Some Constants}\label{sec4}
The different analytical techniques utilized to confirm the irrationality, transcendence, and irrationality measures of many constants are important in the development of other irrationality proofs. Some of these results will be used later on.
\\

\begin{thm} \label{thm4.1} The real numbers \(\pi, \;\zeta (2),\text{ and } \zeta (3)\) are irrational numbers.
\end{thm}

The various irrationality proofs of these numbers are widely available in the open literature. These technique are valuable tools in the theory of irrational numbers, refer to \cite{AR79}, \cite{BF79}, \cite{HD01}, \cite{SJ05}, and others.\\

\begin{thm} \label{thm4.2} For any fixed \(n \in \mathbb{N}\), and the nonprincipal character $\chi \mod 4$, the followings statements are valid.
\begin{enumerate} [font=\normalfont, label=(\roman*)]
\item The real number \(\displaystyle \zeta (2n)=\frac{(-1)^{n+1}2^{2n}B_{2n}}{(2n)!}\pi^{2n}\) is a transcendental number,
\item The real number \(\displaystyle L(2n+1,\chi )=\frac{(-1)^n E_{2n}}{2^{2n+2}(2n)!}\pi^{2n+1}\) is a transcendental number, \\
where \(B_{2n}\text{ and } E_{2n}\) are the Bernoulli and Euler numbers respectively.
\end{enumerate}
\end{thm}

\begin{proof} Apply the Lindemann-Weierstrass theorem to the transcendental number $\pi$. \end{proof}

The first few nonvanishing Bernoulli numbers and Euler numbers are as these.
\begin{enumerate}
\item $B_0=1, B_1=\frac{-1}{2},B_2=\frac{1}{6},B_4=\frac{-1}{30},B_6=\frac{1}{42}, \ldots,$ 
\item $E_0=1, E_2=-1,E_4=5,E_6=-161, \ldots$ .    
\end{enumerate}

The generalization of these results to number fields is discussed in \cite{ZD86}, and related literature.

\begin{thm} \label{thm4.3} {\normalfont (Klinger)} Let $\mathcal{K}$ be a number field extension of degree $k=[\mathcal{K} : \mathbb{Q}$, and discriminant $D=disc (\mathcal{K})$. Then

\begin{enumerate} [font=\normalfont, label=(\roman*)]
\item  If $D>0$, the number field is totally real and $\displaystyle \zeta_{\mathcal{K}}(2n)=r_k\frac{\pi^{2nk}}{\sqrt{D}}$, where $n \geq 1$, and $r_k \in \mathbb{Q}$.
\item  If $D<0$, the number field is totally complex and $\displaystyle \zeta_{\mathcal{K}}(1-2n)=r_k$, where $n \geq 1$, and $r_k \in \mathbb{Q}$.
\end{enumerate}
\end{thm}

\section{Irrationality of the Zeta Constants \(\zeta (2n+1)\)}\label{sec5}
For any integer \(1<s \in \mathbb{N}\), the zeta constant \(\zeta (s)\) is a real number classified as a period since it has a representation as an absolutely convergent integral of a rational function:
\begin{equation}
\zeta (s)=\underset{1>x_1>x_2>\cdots >x_s}{\int }\frac{d x_1}{x_1}\frac{d x_2}{x_2}\cdots \frac{d x_s}{1-x_s} =\sum _{n \geq 1} \frac{1}{n^s},
\end{equation}
where \(s>1\). A few related integral representations are devised in \cite{BF79} to prove the irrationality \(\zeta (2)\) and \(\zeta (3)\). The general idea of a rational or nonrational integral proof of the zeta constant \(\zeta (s)\) for any integer \(s\geq  2\) is probably feasible.\\

\section{The Main Result}\label{sec6}
A different technique using two independent infinite sequences of rational approximations of the two constants \(\zeta _K(s),\text{ and } 1/L(s,\chi )\), which are linearly independent over the rational numbers, will be used to construct an infinite sequence of rational approximations for the zeta constant \(\zeta (2n+1),n \geq 1\). The properties of these sequences, such as sufficiently fast rates of convergence, are then used to derive the
irrationality of any zeta constant \(\zeta (2n+1),n \geq 1\). \\

\begin{proof}  (Proof of Theorem 1.1.) Let \(\chi \neq 1\) be the quadratic character modulo \(q>1,\) and fix an odd integer \(s=2k+1 \geq  3\). By Lemma \ref{lem3.1}, the summatory function of the Dedekind zeta function satisfies the expression
\begin{equation} \label{507}
\zeta (s)\left(\sum _{n \leq  x} \frac{r(n)}{n^s}\right)^{-1}-\frac{1}{L(s,\chi )}=\frac{c_0}{x^{s-1}}+O\left(\frac{1}{x^{s-1/2}}\right)
\end{equation}
for every \(s>1\), and a constant \(c_0>0\). By Theorem \ref{thm4.2}, the real number number \(L(2k+1,\chi )=a \pi ^{2k+1}\), where \(a\in \mathbb{Q}\),
is a transcendental number, so there exists an infinite sequence of rational approximations $\{$\(p_n/q_n:n\in \mathbb{N}\) $\}$ such that
\begin{equation}\label{508}
\left|  \frac{1}{L(s,\chi )}-\frac{p_n}{q_n} \right| <\frac{c_1}{a_nq_n^2} ,
\end{equation}
where $[a_0,a_1,a_2, \ldots ]$ is the continued fraction of $L(s,\chi )^{-1}$, and \(0< c_1<1\) is a constant, see Theorem \ref{thm2.1} or Theorem \ref{thm2.2}. The size of the constant $c_1>0$ in rational approximations is discussed in \cite[p.\ 28]{WM00}. Combining these data, taking absolute value, and using the triangle inequality 
\begin{equation} \label{509}
\left| \left|  x-y \right| \right| \geq \left| \left|  x \right| \right| -\left| \left|  y \right| \right|, 
\end{equation}
lead to the followings.
\begin{eqnarray}\label{510}
\left|  \frac{c_0}{x^{s-1}}+O\left(\frac{1}{x^{s-1/2}}\right) \right| &=&\left|  \zeta (s)\left(\sum _{n \leq  x} \frac{r(n)}{n^s}\right)^{-1}-\frac{1}{L(s,\chi )} \right|  \nonumber  \\                                                                
 &=& \left|  \zeta (s)\left(\sum _{n \leq  x} \frac{r(n)}{n^s}\right)^{-1}-\frac{1}{L(s,\chi)}+\frac{p_n}{q_n}-\frac{p_n}{q_n} \right| \nonumber\\       
&\geq &\left|  \zeta (s)\left(\sum _{n \leq  x} \frac{r(n)}{n^s}\right)^{-1}-\frac{p_n}{q_n}
\right| -\left|  \frac{1}{L(s,\chi )}-\frac{p_n}{q_n} \right| \nonumber\\
&>&\left|  \zeta (s)\left(\sum _{n \leq  x} \frac{r(n)}{n^s}\right)^{-1}-\frac{p_n}{q_n}
\right| -\frac{c_1}{a_nq_n^2} .
\end{eqnarray}

Rewrite it as
\begin{eqnarray} \label{511}
0&<&\left|  \zeta (s)-\frac{p_n}{q_n}\sum _{n \leq  x} \frac{r(n)}{n^s} \right| \\
&<&\left(\frac{c_1}{q_n^2}
+\left|  \frac{c_0}{x^{s-1}}+O\left(\frac{1}{x^{s-1/2}}\right) \right| \right)\sum _{n \leq  x} \frac{r(n)}{n^s} \nonumber .           
\end{eqnarray}
Now, taking \(x \geq  1\) to infinity yields

\begin{equation}\label{512}
0<\left|  \zeta (s)-\frac{p_n}{q_n}\sum _{n \geq 1} \frac{r(n)}{n^s} \right| <\frac{c_2}{a_nq_n^2} ,
\end{equation}
where \(0< c_2<1\) is a constant, and \(s>2\).\\

The rest of the proof is broken up into two cases.\\

\textbf{Case I.} \textit{ Assume that the constant \(\sum _{n=1}^{\infty } r(n)n^{-s}\in \mathbb{Q}\) is a rational number, and that the constant \(\zeta (2k+1)=A/B\) is a rational number.} In this case Theorem \ref{thm4.2} and (\ref{400}) leads to a contradiction.\\

\textbf{Case II.} \textit{ Assume the constant \(\sum _{n=1}^{\infty } r(n)n^{-s}\in \mathbb{R}-\mathbb{Q}\) is an irrational number, and that the constant \(\zeta (2k+1)=A/B\) is a rational number.}\\

By Theorem \ref{thm2.2}, there exists an infinite sequence of rational approximations $\{u_m/v_m:m\in \mathbb{N}\}$ such that
\begin{equation} \label{514}
0<\left|  \sum _{n \geq 1} \frac{r(n)}{n^s} -\frac{u_m}{v_m}\right| <\frac{c_3}{b_mv_m^2} ,
\end{equation}
where $[b_0,b_1,,b_2, \ldots ]$ is the continued fraction of $ \sum _{n \geq 1} r(n)n^{-s}$, and $0< c_3<1$ is a constant. This inequality is equivalent to
\begin{equation}\label{515}
\frac{u_m}{v_m}-\frac{c_3}{b_mv_m^2}<\sum _{n \geq 1} \frac{r(n)}{n^s} <\frac{u_m}{v_m}+\frac{c_3}{b_mv_m^2}.
\end{equation}
Replacing this approximation into inequality (\ref{512}) returns
\begin{equation} \label{516}
0<\left|  \zeta (s)-\frac{p_nu_m}{q_nv_m} \right| \leq\frac{c_4p_n}{b_mv_m^2q_n} +\frac{c_2}{a_nq_n^2} \nonumber ,
\end{equation}
where \(c_4p_n<p_n\) and $0<c_4 <1$ is a constant. By hypothesis, \(\zeta (2k+1)=A/B\) is a rational number, thus, there exists a constant \(c_5\geq 1/B\) such that
\begin{equation} \label{592}
\frac{c_5}{q_nv_m} \leq \left|  \zeta (s)-\frac{p_nu_m}{ q_nv_m} \right| \leq\frac{c_4p_n}{b_mq_nv_m^2} +\frac{c_2}{a_nq_n^2} ,
\end{equation}
this follows from Lemma \ref{lem2.1}. Next, an infinite subsequence of rational approximations 
\begin{equation} \label{519}
\left \{\frac{p_{n_d}u_{m_d}}{q_{n_d}v_{m_d}}: d\geq 1 \right \} \subset \left \{\frac{p_{n}u_{m}}{q_{n}v_{m}}: n,m\geq 1 \right \}
\end{equation}
is generated by the following algorithm. \\

For each $d\geq 1$, and $\varepsilon>0$ is a small number.
\begin{enumerate}
\item If $a_{n_d}> b_{m_d}$, fix the convergent $p_{n_d}/q_{n_d}$ and choose a convergent $u_{m_d}/v_{m_d}$ such that 
\begin{equation} \label{566}
c_6a_{n_d}^{1-\varepsilon}q_{n_d}\leq 2^{w_d}  \leq v_{m_d}\leq 2^{w_d+1}  \leq c_7 a_{n_d}^{1-\varepsilon}q_{n_d}  ,
\end{equation}
or
\item If $a_{n_d}< b_{m_d}$, fix the convergent $p_{n_d}/q_{n_d}$ and and choose a convergent $u_{m_d}/v_{m_d}$ such that 
\begin{equation} \label{567}
c_6b_{m_d}^{1-\varepsilon}v_{m_d}\leq 2^{w_d}   \leq q_{n_d} \leq  2^{w_d+1} \leq c_7 b_{m_d}^{1-\epsilon}v_{m_d}  ,
\end{equation}
where $0<c_6,c_7\leq 2$ are constants. 
\end{enumerate}
 
Replacing the appropiate convergent (\ref{566}) or (\ref{567}) into (\ref{592}) yields
\begin{equation} \label{545}
\frac{c_8}{a_{n_d}^{1-\epsilon}q_{n_d}^2} \leq \left|  \zeta (s)-\frac{p_{n_d}u_{m_d} }{ q_{n_d} v_{m_d}} \right| \leq\frac{c_9p_{n_d}}{q_{n_d}}\frac{1}{a_{n_d}^{2(1-\epsilon)}b_{m_d}q_{n_d}^2 } +\frac{c_2}{a_{n_d}q_{n_d}^2} ,
\end{equation}
or
\begin{equation} \label{546}
\frac{c_8}{b_{m_d}^{1-\epsilon}v_{m_d}^2} \leq \left|  \zeta (s)-\frac{p_{n_d}u_{m_d} }{ q_{n_d} v_{m_d}} \right| \leq\frac{c_9p_{n_d}}{q_{n_d}}\frac{1}{b_{m_d}v_{m_d}^2 } +\frac{c_2}{a_{n_d}b_{m_d}^{2(1-\epsilon)}v_{m_d}^2} ,
\end{equation}
where $0<c_8,c_9,c_{10}<1$ are constants. 

Since the terms
\begin{equation} \label{547}
\frac{c_8p_{n_d}}{q_{n_d}}\frac{1}{a_{n_d}^{2(1-\epsilon)}q_{n_d}^2 } \leq \frac{c_{10}}{a_{n_d}^{2(1-\epsilon)}q_{n_d}^2 }  ,
\end{equation}
and
\begin{equation} \label{548}
\frac{c_8p_{n_d}}{q_{n_d}}\frac{1}{b_{m_d}v_{m_d}^2 } \leq \frac{c_{10}}{b_{m_d}v_{m_d}^2 }  ,
\end{equation}
for all large integers $n_d \geq 1$, see (\ref{508}) and (\ref{514}) respectively, the relation (\ref{545}) or (\ref{546})  reduces to
\begin{eqnarray} \label{584}
\frac{c_8}{a_{n_d}^{1-\epsilon}q_{n_d}^2} &\leq&   \frac{c_9{n_d}}{q_{n_d}}\frac{1}{a_{n_d}^{2(1-\epsilon)}q_{n_d}^2 } +\frac{c_2}{a_{n_d}q_{n_d}^2} \nonumber \\
&\leq & \frac{c_{10}}{a_{n_d}^{2(1-\epsilon)}q_{n_d}^2 }+\frac{c_2}{a_{n_d}q_{n_d}^2},
\end{eqnarray}
or 
\begin{eqnarray} \label{585}
\frac{c_8}{b_{m_d}^{1-\epsilon}v_{m_d}^2} &\leq&\frac{c_{9}p_{n_d}}{q_{n_d}}\frac{1}{b_{m_d}v_{m_d}^2 } +\frac{c_2}{b_{m_d}v_{m_d}^2} \nonumber \\
&\leq & \frac{c_{10}}{b_{m_d}v_{m_d}^2 }+\frac{c_2}{b_{m_d}^{2(1-\epsilon)}v_{m_d}^2}.
\end{eqnarray}

Clearly, the inequality (\ref{584}) or (\ref{585}) leads to a contradiction for infinitely many arge pairs of convergents  \(p_{n_d}/q_{n_d}\) and \(u_{m_d}/v_{m_d}\) as $ n_d, m_d \to \infty$. Ergo, the constant \(\zeta (2k+1)\) is not a rational number. 
\end{proof}

The infinite subsequence of rational approximations
\begin{equation} \label{519}
\frac{p_{n_d}u_{m_d}}{q_{n_d}v_{m_d}} \quad \longrightarrow \quad \zeta (s) \qquad \text{ as }n_d,m_d\longrightarrow \infty ,
\end{equation}
is probably sparse depending on the magnitute of quotiens. However, iterative process in algorithm (\ref{566}) shows that it is infinite, and it suffices to prove the irrationality of $\zeta(2k+1), k \geq 1.$


\section{Irrationality of the Beta Constants \(L(2n,\chi )\)}\label{sec7}
For \(q = 4\) the quadratic symbol is defined by \(\chi (n)=(-1)^{(n-1)/2}\) if \(n \in \mathbb{N}\) is odd, else \(\chi (n)=0\). The corresponding Dedekind zeta function is given by
\begin{equation} \label{620}
\zeta _K(s)=\zeta (s)L(s,\chi )=\frac{1}{4}\sum _{n \geq 1} \frac{r(n)}{n^s} ,
\end{equation}
where the counting function \(r(n)=4\sum _{d | n} \chi (d)=\#\left\{ (a,b):n=a^2+b^2 \right\}\geq 0\) tallies the number of representations of an integer \(n \geq 1\) as sums of two squares, and \(s \in \mathbb{C}\) is a complex number. This is the zeta function of the Gaussian quadratic numbers field \(\mathbb{Q}\left(\sqrt{-1}\right)\). The corresponding $L$-series is
\begin{equation}\label{621}
L(s,\chi )=\sum _{n \geq 1} \frac{\chi (n)}{n^s} =1-\frac{1}{3^s}+\frac{1}{5^s}-\frac{1}{7^s}+\cdots  .
\end{equation}
The evaluation at \(s = 2\) is known as Catalan constant 
\begin{equation}\label{622}
L(2,\chi )=1-\frac{1}{3^2}+\frac{1}{5^2}-\frac{1}{7^2}+\cdots =.915965594177\ldots .
\end{equation} \\

\begin{proof} (Proof of Theorem \ref{thm1.2}.) Use the symmetry of the factorization of the Dedekind zeta function \(\zeta _K(s)=\zeta (s)L(s,\chi )\) with respect
to the zeta function \(\zeta (s)\) and the $L$-function \(L(s,\chi )\) to arrive at the asymptotic formula
\begin{equation} \label{623}
L(s,\chi )\left(\sum _{n \leq  x} \frac{r(n)}{n^s}\right)^{-1}-\frac{1}{\zeta (s)}=\frac{c_0}{x^{s-1}}+O\left(\frac{1}{x^{s-1/2}}\right)
\end{equation}
compare this to (\ref{507}). Now, proceed as before in the proof of Theorem \ref{thm1.1} for the verification of the irrationality of the zeta constant $\zeta (2k+1)$, mutatis mutandis.  
\end{proof}

\section{The \textit{w}-Transform}
A transform as the Laplace transform, Fourier transform, Mellin transform, finite Fourier transform, z-transform, and other related functionals, performs a change of domains to solve certain problems by simpler methods. Likewise, the \textit{\textit{w}-transform} converts some apparently intractable decision problems in the real domain $\R$ to simpler decision problems in the binary domain $\F_2=\{0,1\}$. 

\begin{dfn} 
{\normalfont Let $\alpha \in \R$ be a real number. The \textit{w}-transform is a map $\mathcal{W}: \R \longrightarrow \F_2=\{0,1\}$ defined by
\begin{equation}\label{eq388.44}
\mathcal{W}(\alpha)=
 \lim_{x \to \infty}\frac{1}{2x} \sum_{-x \leq n \leq x}e^{i\alpha n}.
\end{equation}
}
\end{dfn}
The normalization is intrinsic to the number $\pi$. But, it can be modified as needed. The \textit{w}-transform is a point map or equivalently a class map, and it is not invertible. But, inversion is not required in applications to decision problems. 
\begin{lem}\label{lem388.06} For any real number $\alpha \in \R$, the \textit{w}-transform satisfies the followings.

\begin{equation}\label{eq757.37}
\mathcal{W}(2 \pi m\alpha)=
\begin{cases}
1 & \text{ if and only if } \alpha \in \Q,\\
0 & \text{ if and only if } \alpha \notin \Q,
\end{cases}
\end{equation}
for some $m \in \Z$.
\end{lem}

\begin{proof} Given any rational number $\alpha \in \Q$, there is an integer $m \in \Z$ such that $\alpha m \in \Z$, so by definition
\begin{equation}\label{eq388.47}
\mathcal{W}(2 \pi m\alpha)=\lim_{x \to \infty}\frac{1}{2x} \sum_{-x\leq n \leq x}e^{i2  \pi \alpha mn}
=\lim_{x \to \infty}\frac{1}{2x} \sum_{-x\leq n \leq x}1=1.
\end{equation}
The above proves that for some integer $m$, the  sequence 
\begin{equation}
\{2  \pi \alpha mn: n \in \Z\}
\end{equation} 
is not uniformly distributed. While for any irrational number$\alpha \notin \Q$, and any integer $m\ne 0$, the sequence
\begin{equation}
\{2  \pi \alpha mn: n \in \Z\},
\end{equation}
 is uniformly distributed, the proof is the same as the Weil criteria, see \cite[Theorem 2.1]{KN74}.
\end{proof} 
As it is evident, the class function $\mathcal{W}$ maps the class of rational numbers $\Q$ to $1$ and the class of irrational numbers $\I=\R-\Q$ to $0$. The \textit{w}-transform induces an equivalence relation on the set of real numbers $\R=\Q- \I$:
\begin{itemize}
 \item A pair of real numbers $a$ and $b$ are equivalent $a \sim b$ if and only if $\mathcal{W}(2 \pi a)=\mathcal{W}(2 \pi b)$.
 \item A pair of real numbers $a$ and $b$ are not equivalent $a \not \sim b$ if and only if $\mathcal{W}(2 \pi a)\ne\mathcal{W}(2 \pi b)$.
\end{itemize}
The next result takes this concept of class function a little further ahead.

\begin{lem}\label{lem388.99} For any real number $\alpha \in \R$, there is a map $\mathcal{T}: \R \longrightarrow \F_3=\{-1,0,1\}$ that satisfies the followings.

\begin{equation}\label{eq757.37}
\mathcal{T}(2 \pi m\alpha)=
\begin{cases}
1 & \mbox{ if and only if $\alpha$ is rational},\\
0 & \mbox{ if and only if $\alpha$ is algebraic and irrational},\\
-1 &\mbox{ if and only if $\alpha$ is nonalgebraic and irrational},
\end{cases}
\end{equation}
for some $m \in \Z$.
\end{lem}

\begin{lem}\label{lem388.32} For any real number $t \ne k \pi$, $k \in \Z$, and a large integer $x \geq 1$, the finite sum
\begin{enumerate} [font=\normalfont, label=(\roman*)]
 \item $$\sum_{-x\leq n \leq x}e^{i2 tn}=\frac{\sin((2x+1)t)}{\sin(t)}. $$
\item  $$\left |\sum_{-x\leq n \leq x}e^{i2 tn} \right | \leq\frac{1}{| \sin(t) |}. $$
\end{enumerate}
\end{lem}
\begin{proof} (i) Expand the complex exponential sum into two subsums:
\begin{equation}\label{eq957.50}
\sum_{-x\leq n \leq x}e^{i2 tn}=e^{-i2  t}\sum_{0\leq n \leq x-1}e^{-i2  tn}+\sum_{0\leq n \leq x}e^{i2 tn}.
\end{equation}
Lastly, use the geometric series to determine the closed form.	
\end{proof}

 \begin{exa} \label{exa388.37} {\normalfont This demonstration of the \textit{w}-transform verifies that the number $\ln \pi$ is irrational. To confirm this, assume it is a rational number $\ln \pi = r \in \Q^{\times}$, and consider the equivalent equation 
 \begin{equation} \label{eq388.55}
2\pi =2e^{r} .
 \end{equation}
By the Hermite-Lindemann theorem, $e^r$ is transcendental, see \cite[Theorem 1.4]{BA75}, \cite[Theorem 1.2]{WM00}. Take the \textit{w}-transform in both sides to obtain 
\begin{equation} \label{eq388.19}
\mathcal{W}(2 \pi)= \mathcal{W}(2e^r ).
\end{equation} 
The \textit{w}-transforms on the left and right sides are evaluated separately.\\

\textbf{Left Side:} Use the identity $e^{i2 \pi }=1$ to evaluate the the left side of equation (\ref{eq388.19}) as
\begin{equation}\label{eq388.14}
\mathcal{W}(2 \pi)=\lim_{x \to \infty}\frac{1}{2x} \sum_{-x\leq n \leq x}e^{i2  \pi n}=\lim_{x \to \infty}\frac{1}{2x} \sum_{-x\leq n \leq x}1=1.
\end{equation}

\textbf{Right Side:} Use $\sin(e^r)\ne 0$ for an irrational number, and Lemma \ref{lem388.32}, to evaluate the right side of equation (\ref{eq388.19}) as  
\begin{eqnarray}\label{eq388.16}
\mathcal{W}(2e^r )&=& \lim_{x \to \infty}\frac{1}{2x} \sum_{-x \leq n \leq x}e^{i2e^r n} \nonumber \\
&\leq &\lim_{x \to \infty}\frac{1}{2x} \frac{1}{\left | \sin\left (e^r \right )\right |} \\
&=&0 \nonumber.
\end{eqnarray}
The evaluations in (\ref{eq388.14}) and (\ref{eq388.16}) of the \textit{w}-transforms contradict equation (\ref{eq388.19}). Specifically,
\begin{equation} \label{eq388.21}
1=\mathcal{W}(2 \pi)\ne \mathcal{W}(2e^r )=0.
\end{equation}
Therefore, the number $\ln \pi \in \R$ is not a rational number. 
}
\end{exa}

Similar analysis works for other numbers such as $e^{\pi}$, but require additional work to prove the irrationality of $\log_2 \pi$, and $2^{\pi}$.

\begin{rem} {\normalfont  The same result as example \ref{exa388.37} can be proved using the fact that the continued fraction of $e^r$ has one or more arithmetic progressions, \cite{RG73}. which is genenrated by the continued fraction of $e=[2;2,1,2,1,1,4,1,1,6,1,1,8, \ldots]$, see see \cite{EL1744},  \cite[Theorem 3.10]{SJ05}. But, the continued fraction of $\pi=[3;7,15,1,292,1,1,1,2,1,3,1,14, \ldots]$ does not have any known arithmeetic progression. }
\end{rem}

\begin{exe} 
{\normalfont Extend the \textit{w}-transform to a three level map $\mathcal{W}: \R \longrightarrow \F_3$ and prove Lemma \ref{lem388.99}, this requires a deeper understanding of rate of convergence in the Weil criteria to distinguish transcendental numbers from other irrationals.
}
\end{exe}

\begin{exe} 
{\normalfont Prove or disprove the existence of a two level map $\mathcal{N}: \R \longrightarrow \F_2$ defined by 
$$\mathcal{N}( \alpha)=
\begin{cases}
1 & \mbox{ if and only if $\alpha \in \R$ is not normal},\\
0 & \mbox{ if and only if $\alpha \in \R$ is normal}
\end{cases}
$$

}
\end{exe}

\section{Formulas For Zeta Numbers} \label{s190}
For $s \geq 2$, the zeta constant is defined by the Dirichlet series
\begin{equation}
 \zeta(s)=\sum_{n \geq1} \frac{1}{n^s} .
\end{equation}

\begin{lem}  {\normalfont (Euler)} A zeta constant at the even integer argument has an exact Euler formula
\begin{equation} \label{eq257.09}
 \zeta(2n)=(-1)^{n+1}\frac{(2 \pi)^{2n} B_{2n}}{2(2n)!} 
\end{equation} 
in terms of the Bernoulli numbers $B_{2n}$, for $n \geq 1$. 
\end{lem}
\begin{proof} Let $s=2n$, and $B_{2n}(\{x\})$ be the $2n$th Bernouli polynomial, and the corresponding Fourier series
\begin{equation}
\sum_{m \geq 1}\frac{\cos(2m x)}{m^{2n}}=\frac{(-1)^{n+1}(2 \pi)^{2n}}{2} \frac{B_{2n}(\{x\})}{(2n)!}.
\end{equation}
Evaluating at $x=0$ yields $B_{2n}(\{x\})=B_{2n}(0)=B_{2n}.$
\end{proof}
Standard references are \cite[p.\ 18]{CO07}. This formula expresses each zeta constant $\zeta(2n)$ as a rational multiple of $\pi^{2n}$. The formula for the evaluation of the first even zeta constant $\zeta(2)$, known as the Basel problem, was proved by Euler, later it was generalized to all the even integer arguments. Today, there are dozens of proofs, see \cite{CR99}, and \cite[Chapter 6]{SJ03} for an elementary introduction. The first few are 
\begin{multicols}{3}
\begin{enumerate} 
\item $ \displaystyle \zeta(2)=\frac{ \pi^2}{6}, $
\item $ \displaystyle \zeta(4)=\frac{ \pi^4}{90}, $
\item $ \displaystyle \zeta(6)=\frac{ \pi^6}{945} ,$
\end{enumerate}
\end{multicols}

et cetera.\\

In contrast, the evaluation of a zeta constant at an odd integer argument has one or two complicated transcendental power series. A formula for $\zeta(2n+1)$ expresses this constant as a sum of a rational multiple of $\pi^{2n+1}$ and a power series. The derivations involve the Ramanujan series for the zeta function, and appear in \cite[Theorem 1]{GE72}, \cite{GE72B}, \cite{VL06}, \cite{BS17}, et alii. 
The general forms of these formulas are
\begin{equation} 
\zeta(s)=
\begin{cases}
\displaystyle a_n \pi^{4n-1}  -b_n \sum_{n \geq 1} \frac{1}{n^{4n-1}( e^{2 \pi n}-1)}&\text{if $s=4n-1$},\\
\displaystyle a_n \pi^{4n-3}  -b_n \sum_{n \geq 1} \frac{1}{n^{4n-3}( e^{2 \pi n}-1)}-c_n \sum_{n \geq 1} \frac{1}{n^{4n-3}( e^{2 \pi n}+1)}&\text{if $s=4n-3$},
\end{cases}
\end{equation}
where $a_n, b_n, c_n \in \Q$ are rational numbers. 

The precise structure of the rational factor $r_n$ has been known for quite sometimes:
\begin{equation} 
r_n=
\begin{cases}
\displaystyle \frac{2^{2n+1}}{2n(2n+2)} \sum_{0 \leq v \leq n/2}(-1)^v (2n+2-4v) \binom{2n+2}{4v} B_{2v}B_{2n+2-2v}&\text{if $2n+1 \equiv 1 \bmod 4$},\\
\displaystyle \frac{2^{2n}}{(2n+2)!} \sum_{0 \leq v \leq n/2}(-1)^v \binom{2n+2}{2v} B_{2v}B_{2n+2-2v}&\text{if $2n+1 \equiv 3 \bmod 4$}.
\end{cases}
\end{equation}
But, the precise structure of the real number $u_n$ is much more complex, and involves modular forms such as
\begin{equation}
F_{s}(\tau)= \sum_{n \geq 0} \sigma_{-s}(n) q^{n}
\end{equation}
where $\sigma_{-s}(n)=\sum_{d \mid n}d^{-s}$, and $q=e^{i 2 \pi \tau}$. The analysis of the real number $r_n$, for any integer $s \in \Z$, is discussed in \cite[Theorem 1]{GE72}, \cite{GE72B}, \cite{BS17}, etc.

The first few are 
\begin{enumerate} 
\item $ \displaystyle \zeta(3)=\frac{7 \pi^3}{180}  -2 \sum_{n \geq 1} \frac{1}{n^3( e^{2 \pi n}-1)} $,
\item $ \displaystyle \zeta(5)=\frac{ \pi^5}{294}  -\frac{72}{35} \sum_{n \geq 1} \frac{1}{n^5( e^{2 \pi n}-1)}-\frac{2}{35} \sum_{n \geq 1} \frac{1}{n^5( e^{2 \pi n}+1)}, $
\item $ \displaystyle \zeta(7)=\frac{19 \pi^7}{56700}  -2 \sum_{n \geq 1} \frac{1}{n^7( e^{2 \pi n}-1)} ,$
\end{enumerate}
et cetera. These formulas express each zeta constat $\zeta(2n+1)$ as a nearly rational multiple of $\pi^{2n+1}$. These analysis are summarized in a compact formula.

\begin{dfn} \label{dfn257.37}  { \normalfont Let $s \geq 2$ be an integer. The $\pi$-representation of the zeta constant $\zeta(s)=\sum_{n \geq 1}n^{-s}$ is defined by the formula
\begin{equation} 
\zeta(s)=
\begin{cases}
\displaystyle r_n \pi^{s}&\text{if $s=4n, 4n+2$},\\
\displaystyle r_n \pi^{s}  -u_n&\text{if $s=4n-1,4n-3$},
\end{cases}
\end{equation}
where $ r_n \in \Q$ is a rational number and $u_n \in \R$ is a real number.
}
\end{dfn}

\section{Irrational Zeta Numbers} \label{s192}
The zeta function is defined by the series $\zeta(s)= \sum_{n \geq 1}n^{-s}$. The even zeta constant is of the form $\zeta(2n)=r\pi^{2n}$, where $r \ne 0$ is a rational number, see Section \ref{s257} for the actual description. The number $\pi^{2n}$ is irrational for any $n \geq 1$. The irrationality proof for $n=1$ uses the continued fraction of the tangent function $\tan(x)$, the fact that the numbers $\tan(r)$ are irrationals for any nonzero rational number $r \in \Q^{\times}$, and the value $\arctan(1)=\pi/4$ to indirectly show that the continued fraction 
\begin{equation}
\pi=[3;7,15,1,292,1,1,1,2,1,3,1,14, \ldots]
\end{equation}  
is infinite, see \cite[p.\ 129]{BB04}, \cite{LM97}, \cite{NI47}. Later, simpler versions and new proofs were found by several authors, \cite{NI47}, \cite[p.\ 35]{AZ98}, \cite{SJ05}. \\

The irrationality of the first odd zeta constant $\zeta(3)$ was proved by Apery, see, \cite{AR79}. The irrationality of the other odd zeta constants $\zeta(s)$ remain unknown for $s \geq 5$, see \cite{FZ18} and \cite{ZW18}. The $\pi$-representation in Definition \ref{dfn257.37} offers a recursive method for proving the irrationality of $\zeta(2n+1)$ from the known irrationality of $\pi^{2n+1}$ for $n \geq 1$. For example, the irrationality of 
\begin{equation} \label{eq192.59}
\zeta(3), \quad \zeta(5), \quad \zeta(7), \quad \zeta(9), \quad \ldots, 
\end{equation} 
can be derived from the known irrationality of the numbers $\pi^{3}, \pi^{5}, \pi^{7}, \pi^{9}, \ldots $. More generally, this idea can be used to recursively prove the irrationality of $\zeta(s)$ from the known irrationality of $\pi^{s}$ for any integer $s \geq 2$. The inner working of this technique, which requires minimal mathematical knowledge,is demonstrated here for $s=5$. 
\begin{thm} \label{thm192.26}   The number $\zeta(5)$ is irrational.
\end{thm}

 \begin{proof} Suppose that the numbers $1, \pi^5$ and $\zeta(5)$ are linearly dependent over the rational numbers, and consider the equation 
 \begin{equation} \label{eq192.55}
1\cdot a+\pi^5 \cdot b+ \zeta(5) \cdot c=0, 
 \end{equation}
where $(a,b,c)\ne(0,0,0)$ is a nontrivial rational solution. Multiply by $2 \pi $ and the lowest common multiple across the board, and rewrite it in the equivalent form 
\begin{equation} \label{eq192.56}
2 \pi A=-2\left ( B\pi^6+ C\zeta(5) \pi \right )  ,
\end{equation} 
where $A, B, C \in \Z^{\times}$ are integers. To prove the existence or nonexistence of any rational solutions for equation (\ref{eq192.55}), take the \textit{w}-transform in both sides to obtain
\begin{equation} \label{eq192.12}
\mathcal{W}(2 \pi A)= \mathcal{W}\left (-2(B\pi^6 +C\zeta(5) \pi ) \right ).
\end{equation}
The left side and the right side are evaluated separately.\\

\textbf{Left Side.} The evaluation is based on the identity $e^{i2 \pi A}=1$, where $A$ is a fixed integer. The left side evaluation is
\begin{equation}\label{eq192.14}
\mathcal{W}(2 \pi A)=\lim_{x \to \infty}\frac{1}{2x} \sum_{-x\leq n \leq x}e^{i2  \pi An}=\lim_{x \to \infty}\frac{1}{2x} \sum_{-x\leq n \leq x}1=1.
\end{equation}

\textbf{Right Side.} The evalution splits into two cases depending on the values $\sin\left  (B\pi^6 +C\zeta(5) \pi  \right) \ne 0$ or $\sin\left  (B\pi^6 +C\zeta(5) \pi  \right) = 0$.\\

\textbf{Case 1.} $\sin\left  (B\pi^6 +C\zeta(5) \pi  \right) \ne 0$. In this case, Lemma \ref{lem388.32} is applicable. The evaluation for the right side of equation (\ref{eq192.12}) is

\begin{eqnarray}\label{eq192.16}
\mathcal{W}\left (-2(B\pi^6 +C\zeta(5) \pi ) \right )&= &\lim_{x \to \infty}\frac{1}{x} \sum_{-x \leq n \leq x}e^{i2\left (B\pi^6 +C\zeta(5) \pi \right ) n}\\
&\leq &\lim_{x \to \infty}\frac{1}{x}\frac{1}{\sin\left (B\pi^6 +C\zeta(5) \pi \right )} \nonumber\\
&=&0 \nonumber.
\end{eqnarray}
Clearly, these distinct evaluations in equation (\ref{eq192.14}) and in equation (\ref{eq192.16}), that is,
\begin{equation}\label{eq192.49}
1=\mathcal{W}(2 \pi A)\ne\mathcal{W}\left (-2(B\pi^6 +C\zeta(5) \pi ) \right )=0
\end{equation}
contradict equation (\ref{eq192.12}). This implies that equation (\ref{eq192.55}) can not have a nontrivial rational solution $(a,b,c)\ne (0,0,0)$. Hence, the number $\zeta(5) \in \R$ is not a rational number. \\
\end{proof}

\textbf{Case 2.} $\sin\left  (B\pi^6 +C\zeta(5) \pi  \right)=0$. This case implies that equation (\ref{eq192.55}) can have a rational solution $(a,b,c)\ne (0,0,0)$. By Lemma \ref{lem192.93},  
\begin{equation}\label{eq192.92}
\zeta(5)= r_0\pi^5+r_1,
\end{equation}
where $r_0=-B/C, r_1=-m/C \in \Q^{\times}$ are rational numbers. This immediately implies that $\zeta(5)$ is an irrational number. \\

\begin{lem} \label{lem192.93}  If $\sin\left  (B\pi^6 +C\zeta(5) \pi  \right)=0$, then $\zeta(5) =r_0\pi^5 +r_1$, where $r_0=-B/C, r_1=m/C \in \Q^{\times}$ are rational numbers, and $B,C \in \Z^{\times}$ and $m \in \Z$ are integers.
\end{lem}

\begin{proof} The sine function satisfies the relation
\begin{equation}
0= \sin\left  (B\pi^6 +C\zeta(5) \pi  \right)
=\cos \left (B\pi^6 \right)\sin \left ( C\zeta(5) \pi  \right)+\cos \left  (C\zeta(5) \pi  \right)\sin \left ( B\pi^6 \right) .
\end{equation}
Hence, $\tan(B\pi^6)=-\tan(C\zeta(5) \pi)$. Since the tangent function is periodic and one-to-one on the interval $(-\pi/2, \pi/2)$, it implies that $B\pi^6=-C\zeta(5) \pi+m \pi$ for some $m \in \Z$. Equivalently $\zeta(5)=r_0\pi^5+r_1$ with $r_0=-B/C, r_1=-m/C \in \Q^{\times}$. 
\end{proof}

\section{Formulas For Beta Numbers} \label{s196}
The beta function is defined by the Dirichlet series
\begin{equation}
\beta(s)=\sum_{n \geq 1} \frac{\chi(n)}{n^s}=\sum_{n \geq 0} \frac{(-1)^{n}}{(2n+1)^s},
\end{equation}
where $\chi(n)$ is the quadratic symbol, and $s \in C$ is a complex number. A beta constant $\beta(s)$ at an odd integer argument $s=2n+1$ has an exact evaluation.

\begin{lem}  {\normalfont (Euler) } A zeta constant at the even integer argument has the Euler formula
\begin{equation} \label{eq457.09}
 \beta(2n+1)=(-1)^{n+1}\frac{( \pi)^{2n+1} E_{2n}}{4^{n+1}(2n)!} 
\end{equation} 
in terms of the Euler numbers $E_{2n}$, for $n \geq 1$. 
\end{lem}
\begin{proof} Let $s=2n+1$, and $B_{2n+1}(\{x\})$ be the $(2n+1)$th Bernoulli polynomial, and the corresponding Fourier series
\begin{equation}
\sum_{m \geq 1}\frac{\sin(2m x)}{m^{2n+1}}=\frac{(-1)^{n+1}(2 \pi)^{2n+1}}{2} \frac{B_{2n+1}(\{x\})}{(2n+1)!}.
\end{equation}
Evaluating at $x=1/4$ yields 
\begin{equation}
\sum_{m \geq 1}\frac{1}{m^{2n+1}}=\frac{(-1)^{n+1}(2 \pi)^{2n+1}}{2} \frac{B_{2n+1}(1/4)}{(2n+1)!},
\end{equation}
where $2n$th Euler number is defined by
\begin{equation}
E_{2n}=\frac{-4^{2n+1} B_{2n+1}(1/4)}{2n+1}.
\end{equation}
\end{proof}

Some of the standard references are \cite[p.\ 18]{CO07}. This formula expresses each beta constant $\beta(2n+1)$ as a rational multiple of $\pi^{2n+1}$, see \cite{BB87} and related references.  The first few are 
\begin{multicols}{3}
\begin{enumerate} 
\item $ \displaystyle \beta(3)=\frac{ \pi^3}{32}, $
\item $ \displaystyle \beta(5)=\frac{  5\pi^5}{1536} ,$
\item $ \displaystyle \beta(7)=\frac{ 61 \pi^7}{184320},$
\end{enumerate}
\end{multicols}

et cetera. In contrast, the evaluation of a beta constant at an even integer argument can involves the zeta function and a power series, and other complicated formulas, \cite{JL17}, \cite{LF12}, \cite{BW13} et cetera. The derivation for one of the simplest of these formulas is given here.

\begin{lem} \label{lem457.97}  If $\chi(n)$ is the quadratic character, and $s =2k \geq 2$ is an even integer, then the Dirichlet series 
\begin{equation} 
\sum_{n \geq 1} \frac{\chi(n)}{n^s}=\left (1-\frac{1}{2^s} \right )  \zeta(s)  -2 \sum_{n \geq 1} \frac{1}{(4n+3)^s}.
\end{equation}
\end{lem}

\begin{proof} The quadratic character satisfies $\chi(2n+1)=(-1)^n$, and $\chi(2n)=0$. Thus, the Dirichlet series decomposes as 
\begin{eqnarray} 
\sum_{n \geq 1} \frac{\chi(n)}{n^s}&=& \sum_{n \geq 0} \frac{1}{(4n+1)^s} -\sum_{n \geq 0} \frac{1}{(4n+3)^s} \\
&=& \left (1-\frac{1}{2^s} \right )  \zeta(s)-2\sum_{n \geq 0} \frac{1}{(4n+3)^s}\nonumber.
\end{eqnarray}
The last line follows from the identity 
\begin{equation} 
\sum_{n \geq 1} \frac{1}{(2n+1)^s}=\zeta(s)  - \sum_{n \geq 1} \frac{1}{(2n)^s}=\left (1-\frac{1}{2^s} \right ) \zeta(s) 
\end{equation} 
for any complex number $s \in \C$.
\end{proof} 

The first few are 
\begin{enumerate} 
\item $ \displaystyle \beta(2)=\frac{\pi^2}{8}   -2 \sum_{n \geq 1} \frac{1}{(4n+3)^2}, $
\item $ \displaystyle \beta(4)=\frac{7\pi^4}{720}   -2 \sum_{n \geq 1} \frac{1}{(4n+3)^4},$
\item $ \displaystyle \beta(6)=\frac{\pi^6}{960}   -2 \sum_{n \geq 1} \frac{1}{(4n+3)^6},$
\end{enumerate}
et cetera. These analysis are summarized in a compact formula.

\begin{dfn} \label{dfn457.37}  { \normalfont Let $s \geq 2$ be an integer. The $\pi$-representation of the beta constant $\beta(s)$ is defined by the formula
\begin{equation} 
\beta(s)=
\begin{cases}
\displaystyle r \pi^{s}&\text{if $s=2n+1$},\\
\displaystyle r \pi^{s}  -u&\text{if $s=2n$},
\end{cases}
\end{equation}
where $ r \in \Q$ is a rational number and $u \in \R$ is a real number.
}
\end{dfn}
A formula for $\beta(2n+1)$ expresses this constant as a rational multiple of $\pi^{2n+1}$. In contrast, a formula for $\beta(2n)$ expresses this constant as a sum of a rational multiple of $\pi^{2n}$ power series.

\section{Irrational Beta Numbers} \label{s198}
The irrationality of any odd beta constant $\beta(2k+1)$ can be proved by the technique of Lambert, see \cite[p.\ 129]{BB04}. But, the irrationality of any even beta constant $\beta(2n)$ remain unknown for $s=2n \geq 2$, see \cite{NY16} and \cite{ZW18}. The $\pi$-representation in Definition \ref{dfn457.37} offers a recursive method for proving the irrationality of $\beta(2n)$ from the known irrationality of $\pi^{2n}$ for $n \geq 1$. For example, the irrationality of 
\begin{equation} \label{eq457.59}
\beta(2), \quad \beta(4), \quad \beta(6), \quad \beta(8), \quad \ldots, 
\end{equation} 
can be derived from the known irrationality of the numbers $\pi^{2}, \pi^{4}, \pi^{6}, \pi^{8}, \ldots $. More generally, this idea can be used to recursively prove the irrationality of $\beta(s)$ from the known irrationality of $\pi^{s}$ for any integer $s \geq 2$. The inner working of this technique, which is demonstrated here for $s=2$, requires minimal mathematical knowledge. 
\begin{thm} \label{thm457.26}   The number $\beta(2)=\sum_{n \geq 0} (-1)^{n}(2n+1)^{-2}$ is irrational.
\end{thm}

 \begin{proof} Suppose that the numbers $1, \pi^2$ and $\beta(2)$ are linearly dependent over the rational numbers, and consider the equation 
 \begin{equation} \label{eq198.55}
1\cdot a+\pi^2 \cdot b+ \beta(2) \cdot c=0 ,
 \end{equation}
where $(a,b,c)\ne(0,0,0)$ is a nontrivial rational solution. Multiply by $2 \pi$ and the lowest common multiple across the board, and rewrite it in the equivalent form 
\begin{equation} \label{eq198.56}
2 \pi A=-2\left ( B\pi^3+ C\beta(2) \pi \right )  ,
\end{equation} 
where $A, B, C \in \Z^{\times}$ are integers. To prove the existence or nonexistence of any rational solutions for equation (\ref{eq198.55}), take the \textit{w}-transform in both sides to obtain
\begin{equation} \label{eq198.12}
\mathcal{W}(2 \pi A)= \mathcal{W}\left (-2(B\pi^3 +C\beta(2) \pi ) \right ).
\end{equation}
The left side and the right side are evaluated separately.\\

\textbf{Left Side.} The evaluation is based on the identity $e^{i2 \pi A}=1$, where $A$ is a fixed integer. The evaluation of the limit is
\begin{equation}\label{eq198.14}
\mathcal{W}(2 \pi A)=\lim_{x \to \infty}\frac{1}{2x} \sum_{-x\leq n \leq x}e^{i2  \pi An}=\lim_{x \to \infty}\frac{1}{2x} \sum_{-x\leq n \leq x}1=1.
\end{equation}

\textbf{Right Side.} The evalution splits into two cases $\sin\left  (B\pi^3 +C\beta(2) \pi  \right) \ne 0$ and $\sin\left  (B\pi^3 +C\beta(2) \pi  \right) = 0$.\\

\textbf{Case 1.} $\sin\left  (B\pi^3 +C\beta(2) \pi  \right) \ne 0$. In this case, Lemma \ref{lem388.32} is applicable. The evaluation of the limit is 
\begin{eqnarray}\label{eq198.16}
\mathcal{W}\left (-2(B\pi^3 +C\beta(2) \pi ) \right )&= &\lim_{x \to \infty}\frac{1}{x} \sum_{-x \leq n \leq x}e^{i2\left (B\pi^3 +C\beta(2) \pi \right ) n}\\
&\leq &\lim_{x \to \infty}\frac{1}{x}\frac{1}{\sin\left (B\pi^3 +C\beta(2) \pi \right )} \nonumber\\
&=&0 \nonumber.
\end{eqnarray}
Clearly, these distinct evaluations
\begin{equation}\label{eq198.49}
1=\mathcal{W}(2 \pi A)\ne \mathcal{W}\left (-2(B\pi^3 +C\beta(2) \pi ) \right )=0
\end{equation}
contradict equation (\ref{eq198.12}). This implies that equation (\ref{eq198.55}) can not have a nontrivial rational solution $(a,b,c)\ne (0,0,0)$. Hence, the number $\beta(2) \in \R$ is not a rational number. \\
\end{proof}

\textbf{Case 2.} $\sin\left  (B\pi^3 +C\beta(2) \pi  \right)=0$. This case implies that equation (\ref{eq198.55}) can have a rational solution $(a,b,c)\ne (0,0,0)$. By Lemma \ref{lem198.93},  
\begin{equation}\label{eq198.92}
\beta(2)= r_0\pi^2+r_1,
\end{equation}
where $r_0=-B/C, r_1=-m/C \in \Q^{\times}$ are rational numbers. This immediately implies that $\beta(2)$ is an irrational number. \\

\begin{lem} \label{lem198.93}  If $\sin\left  (B\pi^3 +C\beta(2) \pi\right)=0$, then $\beta(2) =r_0\pi^2 +r_1$, where $r_0=-B/C, r_1=m/C \in \Q^{\times}$ are rational numbers, and $B,C \in \Z^{\times}$ and $m \in \Z$ are integers.
\end{lem}

\begin{proof} The sine function satisfies the relation
\begin{equation}
0= \sin\left  (B\pi^3 +C\beta(2) \pi  \right)
=\cos \left (B\pi^3 \right)\sin \left ( C\beta(2) \pi  \right)+\cos \left  (C\beta(2) \pi  \right)\sin \left ( B\pi^3 \right) .
\end{equation}
Hence, $\tan(B\pi^3)=-\tan(C\beta(2) \pi)$. Since the tangent function is periodic and one-to-one on the interval $(-\pi/2, \pi/2)$, it implies that $B\pi^3=-C\beta(2) \pi+m \pi$ for some $m \in \Z$. Equivalently $\beta(2)=r_0\pi^2+r_1$ with $r_0=-B/C, r_1=-m/C \in \Q^{\times}$. 
\end{proof}

\section{Formulas For Beta Numbers} \label{s196}
The beta function is defined by the Dirichlet series
\begin{equation}
\beta(s)=\sum_{n \geq 1} \frac{\chi(n)}{n^s}=\sum_{n \geq 0} \frac{(-1)^{n}}{(2n+1)^s},
\end{equation}
where $\chi(n)$ is the quadratic symbol, and $s \in C$ is a complex number. A beta constant $\beta(s)$ at an odd integer argument $s=2n+1$ has an exact evaluation.

\begin{lem}  {\normalfont (Euler) } A zeta constant at the even integer argument has the Euler formula
\begin{equation} \label{eq457.09}
 \beta(2n+1)=(-1)^{n+1}\frac{( \pi)^{2n+1} E_{2n}}{4^{n+1}(2n)!} 
\end{equation} 
in terms of the Euler numbers $E_{2n}$, for $n \geq 1$. 
\end{lem}
\begin{proof} Let $s=2n+1$, and $B_{2n+1}(\{x\})$ be the $(2n+1)$th Bernoulli polynomial, and the corresponding Fourier series
\begin{equation}
\sum_{m \geq 1}\frac{\sin(2m x)}{m^{2n+1}}=\frac{(-1)^{n+1}(2 \pi)^{2n+1}}{2} \frac{B_{2n+1}(\{x\})}{(2n+1)!}.
\end{equation}
Evaluating at $x=1/4$ yields 
\begin{equation}
\sum_{m \geq 1}\frac{1}{m^{2n+1}}=\frac{(-1)^{n+1}(2 \pi)^{2n+1}}{2} \frac{B_{2n+1}(1/4)}{(2n+1)!},
\end{equation}
where $2n$th Euler number is defined by
\begin{equation}
E_{2n}=\frac{-4^{2n+1} B_{2n+1}(1/4)}{2n+1}.
\end{equation}
\end{proof}

Some of the standard references are \cite[p.\ 18]{CO07}. This formula expresses each beta constant $\beta(2n+1)$ as a rational multiple of $\pi^{2n+1}$, see \cite{BB87} and related references.  The first few are 
\begin{multicols}{3}
\begin{enumerate} 
\item $ \displaystyle \beta(3)=\frac{ \pi^3}{32}, $
\item $ \displaystyle \beta(5)=\frac{  5\pi^5}{1536} ,$
\item $ \displaystyle \beta(7)=\frac{ 61 \pi^7}{184320},$
\end{enumerate}
\end{multicols}

et cetera. In contrast, the evaluation of a beta constant at an even integer argument can involves the zeta function and a power series, and other complicated formulas, \cite{JL17}, \cite{LF12}, \cite{BW13} et cetera. The derivation for one of the simplest of these formulas is given here.

\begin{lem} \label{lem457.97}  If $\chi(n)$ is the quadratic character, and $s =2k \geq 2$ is an even integer, then the Dirichlet series 
\begin{equation} 
\sum_{n \geq 1} \frac{\chi(n)}{n^s}=\left (1-\frac{1}{2^s} \right )  \zeta(s)  -2 \sum_{n \geq 1} \frac{1}{(4n+3)^s}.
\end{equation}
\end{lem}

\begin{proof} The quadratic character satisfies $\chi(2n+1)=(-1)^n$, and $\chi(2n)=0$. Thus, the Dirichlet series decomposes as 
\begin{eqnarray} 
\sum_{n \geq 1} \frac{\chi(n)}{n^s}&=& \sum_{n \geq 0} \frac{1}{(4n+1)^s} -\sum_{n \geq 0} \frac{1}{(4n+3)^s} \\
&=& \left (1-\frac{1}{2^s} \right )  \zeta(s)-2\sum_{n \geq 0} \frac{1}{(4n+3)^s}\nonumber.
\end{eqnarray}
The last line follows from the identity 
\begin{equation} 
\sum_{n \geq 1} \frac{1}{(2n+1)^s}=\zeta(s)  - \sum_{n \geq 1} \frac{1}{(2n)^s}=\left (1-\frac{1}{2^s} \right ) \zeta(s) 
\end{equation} 
for any complex number $s \in \C$.
\end{proof} 

The first few are 
\begin{enumerate} 
\item $ \displaystyle \beta(2)=\frac{\pi^2}{8}   -2 \sum_{n \geq 1} \frac{1}{(4n+3)^2}, $
\item $ \displaystyle \beta(4)=\frac{7\pi^4}{720}   -2 \sum_{n \geq 1} \frac{1}{(4n+3)^4},$
\item $ \displaystyle \beta(6)=\frac{\pi^6}{960}   -2 \sum_{n \geq 1} \frac{1}{(4n+3)^6},$
\end{enumerate}
et cetera. These analysis are summarized in a compact formula.

\begin{dfn} \label{dfn457.37}  { \normalfont Let $s \geq 2$ be an integer. The $\pi$-representation of the beta constant $\beta(s)$ is defined by the formula
\begin{equation} 
\beta(s)=
\begin{cases}
\displaystyle r \pi^{s}&\text{if $s=2n+1$},\\
\displaystyle r \pi^{s}  -u&\text{if $s=2n$},
\end{cases}
\end{equation}
where $ r \in \Q$ is a rational number and $u \in \R$ is a real number.
}
\end{dfn}
A formula for $\beta(2n+1)$ expresses this constant as a rational multiple of $\pi^{2n+1}$. In contrast, a formula for $\beta(2n)$ expresses this constant as a sum of a rational multiple of $\pi^{2n}$ power series.

\section{Irrational Beta Numbers} \label{s198}
The irrationality of any odd beta constant $\beta(2k+1)$ can be proved by the technique of Lambert, see \cite[p.\ 129]{BB04}. But, the irrationality of any even beta constant $\beta(2n)$ remain unknown for $s=2n \geq 2$, see \cite{NY16} and \cite{ZW18}. The $\pi$-representation in Definition \ref{dfn457.37} offers a recursive method for proving the irrationality of $\beta(2n)$ from the known irrationality of $\pi^{2n}$ for $n \geq 1$. For example, the irrationality of 
\begin{equation} \label{eq457.59}
\beta(2), \quad \beta(4), \quad \beta(6), \quad \beta(8), \quad \ldots, 
\end{equation} 
can be derived from the known irrationality of the numbers $\pi^{2}, \pi^{4}, \pi^{6}, \pi^{8}, \ldots $. More generally, this idea can be used to recursively prove the irrationality of $\beta(s)$ from the known irrationality of $\pi^{s}$ for any integer $s \geq 2$. The inner working of this technique, which is demonstrated here for $s=2$, requires minimal mathematical knowledge. 
\begin{thm} \label{thm457.26}   The number $\beta(2)=\sum_{n \geq 0} (-1)^{n}(2n+1)^{-2}$ is irrational.
\end{thm}

 \begin{proof} Suppose that the numbers $1, \pi^2$ and $\beta(2)$ are linearly dependent over the rational numbers, and consider the equation 
 \begin{equation} \label{eq198.55}
1\cdot a+\pi^2 \cdot b+ \beta(2) \cdot c=0 ,
 \end{equation}
where $(a,b,c)\ne(0,0,0)$ is a nontrivial rational solution. Multiply by $2 \pi$ and the lowest common multiple across the board, and rewrite it in the equivalent form 
\begin{equation} \label{eq198.56}
2 \pi A=-2\left ( B\pi^3+ C\beta(2) \pi \right )  ,
\end{equation} 
where $A, B, C \in \Z^{\times}$ are integers. To prove the existence or nonexistence of any rational solutions for equation (\ref{eq198.55}), take the \textit{w}-transform in both sides to obtain
\begin{equation} \label{eq198.12}
\mathcal{W}(2 \pi A)= \mathcal{W}\left (-2(B\pi^3 +C\beta(2) \pi ) \right ).
\end{equation}
The left side and the right side are evaluated separately.\\

\textbf{Left Side.} The evaluation is based on the identity $e^{i2 \pi A}=1$, where $A$ is a fixed integer. The evaluation of the limit is
\begin{equation}\label{eq198.14}
\mathcal{W}(2 \pi A)=\lim_{x \to \infty}\frac{1}{2x} \sum_{-x\leq n \leq x}e^{i2  \pi An}=\lim_{x \to \infty}\frac{1}{2x} \sum_{-x\leq n \leq x}1=1.
\end{equation}

\textbf{Right Side.} The evalution splits into two cases $\sin\left  (B\pi^3 +C\beta(2) \pi  \right) \ne 0$ and $\sin\left  (B\pi^3 +C\beta(2) \pi  \right) = 0$.\\

\textbf{Case 1.} $\sin\left  (B\pi^3 +C\beta(2) \pi  \right) \ne 0$. In this case, Lemma \ref{lem388.32} is applicable. The evaluation of the limit is 
\begin{eqnarray}\label{eq198.16}
\mathcal{W}\left (-2(B\pi^3 +C\beta(2) \pi ) \right )&= &\lim_{x \to \infty}\frac{1}{x} \sum_{-x \leq n \leq x}e^{i2\left (B\pi^3 +C\beta(2) \pi \right ) n}\\
&\leq &\lim_{x \to \infty}\frac{1}{x}\frac{1}{\sin\left (B\pi^3 +C\beta(2) \pi \right )} \nonumber\\
&=&0 \nonumber.
\end{eqnarray}
Clearly, these distinct evaluations
\begin{equation}\label{eq198.49}
1=\mathcal{W}(2 \pi A)\ne \mathcal{W}\left (-2(B\pi^3 +C\beta(2) \pi ) \right )=0
\end{equation}
contradict equation (\ref{eq198.12}). This implies that equation (\ref{eq198.55}) can not have a nontrivial rational solution $(a,b,c)\ne (0,0,0)$. Hence, the number $\beta(2) \in \R$ is not a rational number. \\
\end{proof}

\textbf{Case 2.} $\sin\left  (B\pi^3 +C\beta(2) \pi  \right)=0$. This case implies that equation (\ref{eq198.55}) can have a rational solution $(a,b,c)\ne (0,0,0)$. By Lemma \ref{lem198.93},  
\begin{equation}\label{eq198.92}
\beta(2)= r_0\pi^2+r_1,
\end{equation}
where $r_0=-B/C, r_1=-m/C \in \Q^{\times}$ are rational numbers. This immediately implies that $\beta(2)$ is an irrational number. \\

\begin{lem} \label{lem198.93}  If $\sin\left  (B\pi^3 +C\beta(2) \pi\right)=0$, then $\beta(2) =r_0\pi^2 +r_1$, where $r_0=-B/C, r_1=m/C \in \Q^{\times}$ are rational numbers, and $B,C \in \Z^{\times}$ and $m \in \Z$ are integers.
\end{lem}

\begin{proof} The sine function satisfies the relation
\begin{equation}
0= \sin\left  (B\pi^3 +C\beta(2) \pi  \right)
=\cos \left (B\pi^3 \right)\sin \left ( C\beta(2) \pi  \right)+\cos \left  (C\beta(2) \pi  \right)\sin \left ( B\pi^3 \right) .
\end{equation}
Hence, $\tan(B\pi^3)=-\tan(C\beta(2) \pi)$. Since the tangent function is periodic and one-to-one on the interval $(-\pi/2, \pi/2)$, it implies that $B\pi^3=-C\beta(2) \pi+m \pi$ for some $m \in \Z$. Equivalently $\beta(2)=r_0\pi^2+r_1$ with $r_0=-B/C, r_1=-m/C \in \Q^{\times}$. 
\end{proof}

\begin{exa}
{\normalfont Consider $8\beta(2)=\pi^2+v=\alpha+\kappa$, and $\alpha=1/\pi^2$. The minimal polynomial is
\begin{eqnarray}\label{eq957.50}
g_2(x)&=&\left (x-(\alpha+\kappa) \right ) \left (x-\frac{1}{\alpha} \right )\\
&=&\pi^2 x^2-\left (8\pi^2\beta(2)+ 1 \right ) x+8\beta(2) \nonumber.
\end{eqnarray}
Since $g_2(x) \in \Z[\pi^2,v][x]$ is a polynomial with transcendental coefficients, it follows that the numbers $8\beta(2)= \pi^2+v$ and $1/ \pi^2$ are not algebraic irrational. Therefore, both $8\beta(2)= \pi^2+v$ and $1/ \pi^2$ are transcendental numbers.
}
\end{exa}

\section{Linear Independence Over The Rationals} \label{s957}
Let $\alpha_1, \alpha_2, \ldots,\alpha_d \in \R^{\times}$ be irrational numbers. The existence of rational solutions $c_1, c_2, \ldots,c_d \in \Q$ for the linear equation 
\begin{equation} \label{eq388.49}
c_1\alpha_1+c_2 \alpha_2+ \cdots+c_d\alpha_d=0 
\end{equation}
is an important problem in Diophantine analysis, the general form involving systems of linear equations is discussed in \cite[p.\ 6]{LS71}. Given sufficient information on the parameters, the simplest cases for $d=2$ and possibly $d=3$ can be solved.

\begin{thm} \label{thm957.26}   Let $\alpha \ne r \pi$, $r \in \Q^{\times}$, be an irrational number. Then, the followings numbers are linearly independent over the rational numbers. 
\begin{multicols}{2}
\begin{enumerate} [font=\normalfont, label=(\roman*)]
 \item $1$,  $\alpha$,  and  $\pi$,
\item  $1$,  $\alpha^{-1}$,   and  $\pi$.
\end{enumerate}
\end{multicols}
\end{thm}

\begin{proof} (i) Suppose these numbers are linearly dependent over the rational numbers $\Q$, and consider the equation 
 \begin{equation} \label{eq957.55}
1\cdot a+\alpha \cdot b+\pi \cdot c=0, 
 \end{equation}
where $(a,b,c)\ne(0,0,0)$ is a nontrivial rational solution. Multiply by the lowest common multiple across the board, and rewrite it in the equivalent form 
\begin{equation} \label{eq957.12}
2 \pi C=-2(\alpha B+A)  ,
\end{equation} 
where $A, B,C \in \Z^{\times}$ are integers. To prove the nonexistence of any rational solutions for equation (\ref{eq957.55}) take the \textit{w}-transform in both sides to obtain 
\begin{equation} \label{eq957.19}
\mathcal{W}(2 \pi C)= \mathcal{W}(-2(\alpha B+A )).
\end{equation} 
The \textit{w}-transforms on the left and right sides are evaluated separately.\\

\textbf{Left Side:} Use the identity $e^{i2 \pi C}=1$, where $C$ is an integer, to evaluate the the left side of equation (\ref{eq957.19}) as
\begin{equation}\label{eq957.14}
\mathcal{W}(2 \pi C)=\lim_{x \to \infty}\frac{1}{2x} \sum_{-x\leq n \leq x}e^{i2  \pi Cn}=\lim_{x \to \infty}\frac{1}{2x} \sum_{-x\leq n \leq x}1=1.
\end{equation}

\textbf{Right Side:} Use $\sin(\alpha B+A)\ne 0$ for an irrational number, and Lemma \ref{lem388.32}, to evaluate the right side of equation (\ref{eq957.12}) as  
\begin{eqnarray}\label{eq957.16}
\mathcal{W}(-2(\alpha B+A ))&=& \lim_{x \to \infty}\frac{1}{2x} \sum_{-x \leq n \leq x}e^{-i2(\alpha B+A) n} \nonumber \\
&\leq &\lim_{x \to \infty}\frac{1}{2x}\frac{1}{\left |\sin\left (\alpha B+A \right ) \right |}\\
&=&0 \nonumber.
\end{eqnarray}
The evaluations in (\ref{eq957.16}) and (\ref{eq957.14}) of the \textit{w}-transforms 
\begin{equation} \label{eq957.22}
1=\mathcal{W}(2 \pi B)\ne \mathcal{W}(-2(\alpha B+A ))=0
\end{equation}
contradict equation (\ref{eq957.12}). Therefore, equation (\ref{eq957.55}) can not have a nontrivial rational solution $(a,b,c)\ne (0,0,0)$. The proof of (ii) is similar.
\end{proof}

\begin{thm} \label{thm957.46}   Let $\alpha \ne r \pi$, $r \in \Q^{\times}$, be an irrational number. Then, the followings numbers are linearly independent over the rational numbers. 
\begin{multicols}{2}
\begin{enumerate} [font=\normalfont, label=(\roman*)]
 \item  $1$,  $\alpha$,   and  $\pi^{-1}$,
 \item  $1$,  $\alpha^{-1}$,   and  $\pi^{-1}$.
\end{enumerate}
\end{multicols}
\end{thm}

\begin{proof} (i) Suppose these numbers are linearly dependent over the rational numbers $\Q$, and consider the equation 
 \begin{equation} \label{eq957.75}
1\cdot a+\alpha \cdot b+\pi^{-1} \cdot c=0, 
 \end{equation}
where $(a,b,c)\ne(0,0,0)$ is a nontrivial rational solution. Rewrite it in the equivalent form 
\begin{equation} \label{eq957.22}
2 \pi =\frac{-2c}{\alpha b+a}  ,
\end{equation} 
where $a, b,c \in \Q^{\times}$ are rational numbers. To prove the nonexistence of any rational solutions for equation (\ref{eq957.55}) take the \textit{w}-transform in both sides to obtain 
\begin{equation} \label{eq957.39}
\mathcal{W}(2 \pi )= \mathcal{W}\left (\frac{-2c}{\alpha b+a} \right ).
\end{equation} 
The \textit{w}-transforms on the left and right sides are evaluated separately.\\

\textbf{Left Side:} Use the identity $e^{i2 \pi }=1$ to evaluate the the left side of equation (\ref{eq957.39}) as
\begin{equation}\label{eq957.34}
\mathcal{W}(2 \pi )=\lim_{x \to \infty}\frac{1}{2x} \sum_{-x\leq n \leq x}e^{i2  \pi n}=\lim_{x \to \infty}\frac{1}{2x} \sum_{-x\leq n \leq x}1=1.
\end{equation}

\textbf{Right Side:} Use $\sin\left (\frac{-c}{\alpha b+a} \right )\ne 0$ for an irrational number $\alpha$, and Lemma \ref{lem388.32}, to evaluate the right side of equation (\ref{eq957.39}) as  
\begin{eqnarray}\label{eq957.36}
\mathcal{W}\left (\frac{-2c}{\alpha b+a} \right )&=& \lim_{x \to \infty}\frac{1}{2x} \sum_{-x \leq n \leq x}e^{-i2\left (\frac{-2c}{\alpha b+a} \right ) n} \nonumber \\
&\leq &\lim_{x \to \infty}\frac{1}{2x}\frac{1}{\left |\sin\left (\frac{-c}{\alpha b+a} \right ) \right |}\\
&=&0 \nonumber.
\end{eqnarray}
The evaluations in (\ref{eq957.36}) and (\ref{eq957.34}) of the \textit{w}-transforms 
\begin{equation} \label{eq957.22}
1=\mathcal{W}(2 \pi B)\ne \mathcal{W}\left (\frac{-2c}{\alpha b+a} \right )=0
\end{equation}
contradict equation (\ref{eq957.39}). Therefore, equation (\ref{eq957.75}) can not have a nontrivial rational solution $(a,b,c)\ne (0,0,0)$. The proof of (ii) is similar.
\end{proof}

\begin{cor} \label{cor957.86} For any rational number $r \in \Q$, the followings statements are valid. 
\begin{multicols}{2}
\begin{enumerate} [font=\normalfont, label=(\roman*)]
 \item $\alpha \ne r\pi$;
\item $\alpha^{-1}\ne r\pi$,
\item $\alpha\ne r\pi^{-1}$,
 \item $\alpha^{-1}\ne r\pi^{-1}$.
\end{enumerate}
\end{multicols}
\end{cor}
\begin{proof} (i) By Theorem \ref{thm957.26}, the equation $1\cdot a+\alpha \cdot b+\pi \cdot c=0$ has no rational solutions $(a,b,c)=(0, b,c)\ne (0,0,0)$. 
\end{proof}

\section{The Sum $e+\pi$ And Product $e \pi$ } \label{s992}
\begin{lem} \label{lem957.36}   The numbers $e +\pi$ and $e \cdot\pi$ are irrational numbers.
\end{lem}
\begin{proof} (i) By Theorem \ref{thm957.26}, the equation $1\cdot a+e \cdot b+\pi \cdot c=0$ has no rational solutions $(a,b,c)=(0, b,c)\ne (0,0,0)$. Therefore, $e+\pi=r_0$ has no rational solution $r_0 \in \Q$. \\

(ii) By Theorem \ref{thm957.46}, the equation $1\cdot a+e \cdot b+\pi^{-1} \cdot c=0$ has no rational solutions $(a,b,c)=(0, b,c)\ne (0,0,0)$. Therefore, $e=r_1 \pi^{-1}$ has no rational solution $r_1 \in \Q$.
\end{proof}

\begin{lem} \label{lem957.36}   The numbers $e +\pi$ and $e \pi$ are transcendental, (nonalgebraic irrational numbers).
\end{lem}

\begin{proof} (i) The irrational numbers $e+\pi$ and $e^{-1}\pi$ are the unique roots of the polynomial
\begin{eqnarray}\label{eq957.50}
f(x)&=&\left (x-(e+\pi) \right ) \left (x-e^{-1}\pi \right )\\
&=&x^2-(e+\pi+e^{-1} \pi) x+ \pi+e^{-1} \pi^{2}\nonumber\\
&=&ex^2-(e^2+e\pi + \pi) x+ e\pi+ \pi^{2}\nonumber.
\end{eqnarray}
Since $f(x) \in \Z[e,\pi][x]$ is a polynomial with transcendental coefficients, it follows that the numbers $e+ \pi$ and $e^{-1} \pi$ are not algebraic irrational. Therefore, both $e+ \pi$ and $e^{-1} \pi$ are transcendental numbers, (the roots of a nonalgebraic polynomial).

(ii) The irrational numbers $e^{-1}+\pi$ and $e\pi$ are the unique roots of the polynomial
\begin{eqnarray}\label{eq957.50}
g(x)&=&\left (x-(e^{-1}+\pi) \right ) \left (x-e\pi \right )\\
&=&x^2-(e^{-1}+\pi +e \pi) x+ \pi+e^{-1} \pi^{2}\nonumber\\
&=&ex^2-(1+e^2\pi +e\pi) x+ e\pi+ \pi^{2}\nonumber.
\end{eqnarray}
Since $g(x) \in \Z[e,\pi][x]$ is a polynomial with transcendental coefficients, it follows that the numbers $e^{-1}+ \pi$ and $e \pi$ are not algebraic irrational. Therefore, both $e^{-1}+ \pi$ and $e \pi$ are transcendental numbers, (the roots of a nonalgebraic polynomial).
\end{proof}

\section{Sums and Products of Algebraic And Nonalgebraic Numbers} \label{s462}
 The algebraic closure of the rational numbers consists of all solutions of rational polynomials equations. The subset of real numbers is denoted by 
$$\overline{\Q}=\{\alpha \in \R: f(\alpha)=0 \text{  and  } f(x)\in \Q[x]\}.$$

 \begin{dfn} {\normalfont An irrational number $\alpha \in \C^{\times}$ is called \textit{algebraic irrational} if and only if there is a rational polynomial $f(x)\in \Q[x]$ such that $f(\alpha)=0$. Otherwise, it is called \textit{nonalgebraic irrational} or transcendental.}
 \end{dfn}
 
\begin{dfn} {\normalfont The subset of algebraic irrational numbers is defined by $$\A=\{\alpha \in \R:\alpha \text{ is irrational and } f(\alpha)=0 \}$$ for some rational polynomial $f(x)\in \Q[x]$.}
\end{dfn}
The subset of numbers $\A$ is a proper subset of the set of algebraic integers, that is, $\A \subset \overline{\Q}.$
\begin{dfn} {\normalfont The subset of nonalgebraic irrational numbers is defined by $$\T=\{\alpha \in \R:\alpha \text{ is irrational and } f(\alpha)\ne 0 \}$$ for any rational polynomial $f(x)\in \Q[x]$.}
\end{dfn}
 
\begin{thm} \label{thm462.33} The subsets $\A$ and $\T$ have the followings properties.
\begin{enumerate} [font=\normalfont, label=(\roman*)]
 \item The subset $\A$ of algebraic irrational numbers is pseudo ring without rational numbers $\Q$.
\item  The subset $\T$ of nonalgebraic irrational numbers is pseudo ring without algebraic rational numbers $\overline{\Q}$.
\end{enumerate}
\end{thm}

\begin{proof} (ii) Take a pair of nonalgebraic irrational numbers $\alpha, \beta \in \T$ such that $\alpha\beta \notin \overline{\Q}$. Then, by Lemma \ref{lem462.08}, the sum $\alpha+ \beta \in \T$, and the product $\alpha\beta \in \T$ are nonalgebraic irrationals. The condition $\alpha\beta \notin \overline{\Q}$ implies that the subset $\T$ does not contain the algebraic rational numbers $\overline{\Q}$.
\end{proof}

A new subset of numbers, which is a ring without units, was defined in \cite{KZ01}. This subset is a proper subset of the union of algebraic irrationals and nonalgebraic irrationals numbers. 
\begin{equation}
\mathscr{P} =\{periods \}\subset  \A \cup \T.
\end{equation}
Surprisingly, the set of periods $\mathscr{P}$ is a countable set. \\

Let $\alpha,\beta \in \R^{\times}$ be a real numbers. The sum $\alpha+\beta$ is either a rational number, an irrational number or a transcendental number depending on the property of the number $\alpha $ or $\beta$. This simple observation is used below. 

\begin{lem} \label{lem462.08} Let $\alpha \in \R^{\times}$ be a transcendental number, and let $\beta \in \R^{\times}$ be a real number. Then, the number $\alpha+\beta \in \R^{\times}$ is a transcendental number.
\end{lem}

\begin{proof} The real numbers $\alpha+\beta$ and $1/\alpha$ are the unique roots of the polynomial
\begin{eqnarray}\label{eq957.50}
f(x)&=&\left (x-(\alpha+\beta) \right ) \left (x-\frac{1}{\alpha} \right )\\
&=&x^2-\left (\alpha+\beta+\frac{1}{\alpha} \right ) x+\left (\alpha+\beta \right )\frac{1}{\alpha}\nonumber \\
&=&\alpha x^2-\left (\alpha^2+ \alpha \beta+1 \right ) x+\alpha+\beta \nonumber.
\end{eqnarray}
Since $f(x) \in \Z[\alpha][x]$ is a polynomial with transcendental coefficients, it follows that the numbers $\alpha+\beta$ and $1/\alpha$ are not 

algebraic irrational. Therefore, both are transcendental numbers.
\end{proof}

\end{document}